\documentclass[12pt,english]{article}
\usepackage[T1]{fontenc}
\usepackage{amsfonts}
\usepackage{amssymb}
\usepackage{stmaryrd}
\usepackage{amsthm,amsmath}
\usepackage{a4wide}
\usepackage{shadow}
\usepackage{ulem}
\usepackage{enumerate}
\usepackage[dvips]{graphics}
\usepackage{lastpage}
\usepackage[perpage]{footmisc}
\usepackage{color}

\font\tensym=msbm10
\font\sevensym=msbm7
\font\fivesym=msbm5
\newfam\symfam
\textfont\symfam=\tensym
\scriptfont\symfam=\sevensym
\scriptfont\symfam=\fivesym

\begin{document}
\bibliographystyle{plain}
\title{\bf{On the rate of convergence in the central limit theorem for martingale difference sequences}}
\date{24 Mars 2004}
\author{Lahcen OUCHTI\thanks{LMRS, UMR 6085, Universit\'e de Rouen, site Colbert 76821
Mont-Saint-Aignan cedex, France.\qquad\qquad\qquad E-mail: lahcen.ouchti@univ-rouen.fr}}
\maketitle
\newtheorem{theoreme}{Theorem}
\newtheorem{lemme}{Lemma}
\newtheorem{cor}{Corollary}
\newtheorem{exemp}{Exemple}
{\textsc{Abstract:}} We established the rate of convergence in the
central limit theorem for stopped sums of
 a class of martingale difference sequences.
\begin{center}\textrm{\bf{Sur la vitesse de convergence dans le th\'eor\`eme limite
central pour les diff\'erences de martingale}}\end{center}
\par {\textsc{R\'esum\'e:}} On \'etablit la vitesse de convergence
dans le th\'eor\`eme limite central pour les sommes arr\^et\'ees
issues d'une classe de suites de diff\'erences de martingale.
 \vspace{+0.4cm}

\section{ Introduction}
 Let $(X_{i})_{i\in\mathbb{N}}$ be a sequence of random variables defined on a probability space
 $(\Omega,\,\, \mathcal {F},\,\,\mathbb P)$. We shall say that
 $(X_{k})_{k\in\mathbb {N}}$ is a martingale difference sequence if, for any  $k\ge 0$
 \begin{enumerate}
 \item $\mathbb E \{|X_{k}|\}<+\infty $.
 \item $\mathbb E \{X_{k+1}|\mathcal {F}_{k}\}=0$, where $\mathcal {F}_{k}$ is the $\sigma$-algebra generated by
 $X_{i}, i\le k$.
\end{enumerate}
\par For each integer  $n\ge 1$ and  $x$ real number, we denote
 $$S_{0}=0,\quad S_{n}=\sum_{i=1}^{n} X_{i},
\quad \phi(x)=\frac{1}{\sqrt{2\pi}}\int_{-\infty}^{x}\exp(-\frac{t^{2}}{2}) dt,
 \quad \sigma_{n-1}^{2}=\mathbb E \{X_{n}^{2}|\mathcal {F}_{n-1}\},$$
 $$\nu(n)=\inf\{k\in\mathbb{N}^{*}/
\,\,\sum_{i=0}^{k}\sigma_{i}^{2}\ge n\},\quad S_{\nu(n)}^{2}=\sum_{k=1}^{+\infty}S_{k}^{2}\,I_{\nu(n)=k},
\quad \sigma_{\nu(n)}^{2}=\sum_{k=1}^{+\infty}\sigma_{k}^{2}\,I_{\nu(n)=k},$$
$$\quad F_{n}(x)=
\mathbb P(S_{\nu(n)}\le x\sqrt{n}),\quad S'_{\nu(n)}=S_{\nu(n)}+\sqrt{\gamma(n)}X_{\nu(n)+1},
\quad H_{n}(x)=\mathbb P(S'_{\nu(n)}
\le x\sqrt{n}),$$
and $\gamma(n)$ is a random variable such that
\begin{equation}\label{p}
\sum_{i=0}^{\nu(n)-1}\sigma_{i}^{2}+\gamma(n)\sigma_{\nu(n)}^{2}=n\quad\,\,\,p.s.
\end{equation}
If the random variables $X_{i}$ are independent and identically distributed with
$\mathbb E X_{i}=0$ and $\mathbb E X_{i}^{2}=1$,
 we have by the central limit theorem (CLT)
$$\lim_{n\to +\infty} \sup_{x\in\mathbb {R}}|\mathbb P(S_{n}\le x\sqrt{n})-\phi(x)|=0.$$
By the theorem of Berry (\cite{Br}, 1941) and Esseen (\cite{Es},
1942), if moreover, $\mathbb E |X_{i}^{3}|<+\infty$, the rate of
convergence in the limit is of order $n^{-\frac{1}{2}}$. If
$(X_{i})_{i\in\mathbb N}$ is an ergodic martingale difference
sequence with $\mathbb E X_{i}^{2}=1$, by the theorem of
Billingsley (\cite{Bill}, 1968) and Ibragimov ((\cite{Ibr}, 1963),
see also (\cite{Hall}, 1980)) we have the CLT. The rate of
convergence can, however, be arbitrarily slow even if $X_{i}$ are
bounded and $\alpha$-mixing (cf  \cite {Vol}). There are several
results showing that with certain assumption on the conditional
variance $\mathbb E (X_{i}^{2}|\mathcal F_{i-1})$, the rate of
convergence becomes polynomial (Kato (\cite {Ka}, 1979), Grams
(\cite {Grams}, 1972), Nakata (\cite {Nakata}, 1976),
 Bolthausen (\cite{Bol}, 1982),
 Haeusler (\cite {Haeusler}, 1988), \ldots).
\par In 1963, Ibragimov \cite{Ibr} has shown that for $X_{i}$ uniformly bounded, if instead of usual sums $S_{n}$,
the stopped sums $S_{\nu(n)}$ or $S'_{\nu(n)}$ are considered, one gets the rate of convergence of order
$n^{-\frac{1}{4}}$; the only assumption
beside boundednes is that $\sum_{i=0}^{+\infty}\, \sigma_{i}^{2}$ diverge to infinity a.s.
\par In the present paper we give a rate of convergence for a larger class of martingale difference sequences,
the Ibragimov's case will be a particular one.
 \section{Main result}
 \par We consider a sequence $(X_{i})_{i\in\mathbb{N}}$ of square integrable martingale differences.
\begin{theoreme}\label{theorem}
 If the series $\sum_{i=0}^{+\infty}\, \sigma_{i}^{2}$ diverges a.s.  and if there exists a nondecreasing
 sequence $(Y_{i})_{i\in\mathbb N}$ adapted to the filtration
 $(\mathcal F_{i},\,\, i\in\mathbb N)$ such that, for all
$i\in\mathbb{N}^{*}$
$$\mathbb E (Y_{i}^{4})< +\infty,\,\,\,\,1 \le Y_{i}\quad \textrm{and}\quad
\mathbb E (|X_{i}|^{3}|\mathcal F_{i-1})\le Y_{i-1}\, \sigma_{i-1}^{2}\,\,\,\, a.s.$$
then for all n sufficiently large
\begin{equation}\label{1}
\sup_{x\in\mathbb{R}}\biggl|F_{n}(x)-\phi(x)
\biggr| \leq \frac{{a_{n}}^{\frac{1}{2}}}{\pi {n}^{\frac{1}{4}}}\biggl(11+\frac{3}{4n^{\frac{1}{4}}}+\frac{2}{9n^{\frac{1}{2}}}+
\frac{1}{8n^{\frac{3}{4}}}\biggr),
\end{equation}
\begin{equation}\label{2}
\sup_{x\in\mathbb{R}}\biggl|H_{n}(x)-\phi(x)\biggr| \leq \frac{{a_{n}}^{\frac{1}{2}}}{\pi {n}^{\frac{1}{4}}}\biggl
(11+\frac{9}{4n^{\frac{1}{4}}}+\frac{2}{9n^{\frac{1}{2}}}+
\frac{1}{8n^{\frac{3}{4}}}\biggr)
\end{equation}
where $ a_{n}= (\mathbb E Y_{\nu(n)}^{4})^{\frac{1}{2}}$.
\end{theoreme}
If we put $Y_{i}=M$ a.s. where $M>0$ is a constant, one obtains the following corollaries:
\begin{cor}
 If the series $\sum_{i=0}^{+\infty}\, \sigma_{i}^{2}$ diverges a.s. and there exists  $M>0$ such that, for all
 $i\in\mathbb{N}^{*}$,\,\,$\mathbb E (|X_{i}|^{3}|\mathcal F_{i-1})\le M \,\mathbb E (X_{i}^{2}|\mathcal F_{i-1})$ a.s. then there
  is a constant $0<c_{M}<+\infty$
\begin{equation}\label{3}
\sup_{x\in\mathbb{R}}\biggl|F_{n}(x)-\phi(x)
\biggr| \leq \frac{c_{M}}{ n^{\frac{1}{4}}},
\end{equation}
\begin{equation}\label{4}
\sup_{x\in\mathbb{R}}\biggl|H_{n}(x)-\phi(x)\biggr|\leq \frac{c_{M}}{n^{\frac{1}{4}}}.
\end{equation}
\end{cor}
\begin{cor} If there exists $0<\alpha\le M<+\infty$ satisfying $\sigma_{i-1}^{2}\ge \alpha$ and
$\mathbb E (|X_{i}|^{3}|\mathcal F_{i-1})\le M$ a.s. for all $i\in\mathbb{N}^{*}$,
then there is a constant $0<c_{(\alpha,\, M)}<+\infty$  such that $(\ref {3})$ and $(\ref{4})$ hold.
\end{cor}
Moreover, if we suppose that $(X_{i})_{i\in\mathbb {N}}$ is uniformly bounded, we obtain the result of Ibragimov \cite {Ibr}.
\begin{cor} If the series $\sum_{i=0}^{+\infty}\, \sigma_{i}^{2}$ diverges a.s.
and $|X_{i}|\le M<+\infty$ a.s. for all $i\ge 0$, then $(\ref{3})$ and $(\ref{4})$ hold.
\end{cor}

\hspace{-0.6cm}{\textbf{Example.}} Let $A=(A_{k})_{k\in\mathbb {N}}$ be a sequence of real valued random variables such that
$\sup_{k\in\mathbb {N}}\mathbb{E}(A_{k}^{4})^{1/4}=\beta <\infty$
and consider an arbitrary sequence of variables
$\zeta=(\zeta_{k})_{k\in\mathbb {N}^{*}}$ with zero means, unit
variances, bounded third moments and which are also independent
of $A$. We definie $X=(A_{k-1}\zeta_{k})_{k\in\mathbb {N}^{*}}$
and $\mathcal {F}_{k}$ the $\sigma$-algebra generated by
$A_{0},\,A_{1},\dots,\,A_{k}$.
\par Clearly $(X_{k},\,\,\mathcal {F}_{k},\,\,k\in\mathbb{N}^{*})$ is
a martingale difference sequence, and for all $k\in\mathbb {N}^{*}$,
$$\mathbb {E}(A_{k-1}^{2}\zeta_{k}^{2}|\mathcal {F}_{k-1})=A_{k-1}^{2}\qquad \textrm{a.s.,}$$
$$\mathbb {E}(|A_{k-1}\zeta_{k}|^{3}\mathcal {F}_{k-1})\leq |A_{k-1}|
\sup_{i\in{\mathbb {N}}^{*}}\mathbb {E}(|\zeta_{i}|^{3})\,A_{k-1}^{2}\qquad \textrm{a.s..}$$
If $(|A_{k}|)_{k\in\mathbb {N}}$
is nondecreasing, then using Theorem \ref{theorem}, one obtains
$$\sup_{x\in\mathbb {R}}\biggl|F_{n}(x)-\phi(x)\biggr|\leq
c\beta\frac{\sup_{k\in{\mathbb{N}}^{*}}\mathbb
{E}(|\zeta_{k}|^{3})^{\frac{1}{4}}}{n^{\frac{1}{4}}}$$ where $c$
is a positive constant.

\section{Proof of Theorem}
According to Esseen's theorem ( see, e.g., (\cite{B. V.}, 1954) p. 210 and (\cite{Lo}, 1955) p. 285), for all $y>0$,
\begin{equation}\label{e}
\sup_{x\in\mathbb{R}}\biggl|F_{n}(x)- \phi(x)\biggr| \le \frac{1}{\pi}
\int_{-y}^{y}\biggl|\mathbb E \biggl\{\exp(\frac{itS_{\nu(n)}}{\sqrt{n}})\biggr\}-
\exp(-\frac{t^{2}}{2})\biggr|\,\frac{dt}{|t|}\,+ \frac{24}{\pi \sqrt{2\pi}y}.
\end{equation}
\par Below we shall prove the following inequalities
\begin{equation}
\biggl|\mathbb E \biggl\{\exp\biggl(\frac{itS_{\nu(n)}}{\sqrt{n}}+\frac{t^{2}}{2n}
\sum\limits_{p=0}^{\nu(n)-1}\,\sigma_{p}^{2}\biggr)\biggr\}-1\biggr| \leq
 a_{n}e^{\frac{t^{2}}{2}}\biggl(\frac{|t|}{3\sqrt{n}}+\frac{t^{2}}{4n}+\frac{a_{n}|t|^{3}}{3n^{\frac{3}{2}}}
+\frac{a_{n} t^{4}}{4n^{2}}\biggr),
\end{equation}
\begin{equation}
\biggl|\mathbb E \biggl\{\exp\biggl(\frac{itS_{\nu(n)}}{\sqrt{n}}+
\frac{t^{2}}{2n}\sum\limits_{p=0}^{\nu(n)-1}\,\sigma_{p}^{2}\biggr)\!\biggr\}
-\mathbb E \biggl\{\exp\biggl(\frac{itS_{\nu(n)}}{\sqrt{n}}+
 \frac{t^{2}}{2}\!\biggr)\!\biggl\}\!\biggr| \leq
 \frac{a_{n}t^{2}}{2n}\exp(\frac{t^{2}}{2}),
\end{equation}
\begin{equation}
  \biggl|\mathbb E \biggl\{\exp\biggl(\frac{itS_{\nu(n)}}{\sqrt{n}}\biggr)\biggr\}-
\mathbb E \biggl\{\exp\biggl(\frac{itS'_{\nu(n)}}{\sqrt{n}}\biggr)\biggr\}\biggr|\le \frac{3a_{n}t^{2}}{2n}
\end{equation}
where $a_{n}=(\mathbb E Y_{\nu(n)}^{4})^{\frac{1}{2}}$.
\subsection{Proof of the Inequality (7)}
 We have
\begin{align*}
&\mathbb E \biggl\{\exp\biggl(\frac{itS_{\nu(n)}}{\sqrt{n}}+\frac{t^{2}}{2n}
\sum_{p=0}^{\nu(n)-1}\,\sigma_{p}^{2}\biggr)\biggr\}-1\\
&=\sum_{k=1}^{+\infty}\mathbb E \biggl\{\biggl(\exp\biggl(\frac{itS_{k}}{\sqrt{n}}+\frac{t^{2}}{2n}
\sum_{p=0}^{k-1}\,\sigma_{p}^{2}\biggr)-1\biggr)I_{\nu(n)=k}\biggr\}\\
&=\sum_{k=1}^{+\infty}\sum_{j=1}^{k}\mathbb E \biggl\{\exp\biggl(\frac{itS_{j-1}}{\sqrt{n}}+
\frac{t^{2}}{2n}\sum_{p=0}^{j-1}\sigma_{p}^{2}\biggr)
\biggl(e^{\frac{itX_{j}}{\sqrt{n}}}-e^{-\frac{t^{2}\sigma_{j-1}^{2}}{2n}}\biggr)I_{\nu(n)=k}\biggr\}.
\end{align*}
\par For real $x$, put
$$e^{ix}=1+ix+\frac{(ix)^{2}}{2}+u(x),\quad
e^{-x}=1-x+\beta(x)\frac{x^{2}}{2}\qquad\quad(*)$$ 
It is easily seen that, for all $x\in\mathbb{R}$
$$|u(x)|\le \frac{|x|^{3}}{6},\quad |u(x)|\le \frac{x^{2}}{2},\,\, \begin{text}{and}\end{text}
\,\,\,\,\,\,|\beta(|x|)|\le 1.$$
Observing that the random variable
$W_{j-1}^{n}=\exp\biggl(\frac{itS_{j-1}}{\sqrt{n}}+
\frac{t^{2}}{2n}\sum\limits_{p=0}^{j-1}\sigma_{p}^{2}\biggr)$ is
measurable with respect to the $\sigma$-algebra $\mathcal
{F}_{j-1}$ and using the identities ($*$), we obtain
\begin{align}\label{dalib}
&\mathbb E \biggl\{\exp\biggl(\frac{itS_{\nu(n)}}{\sqrt{n}}+\frac{t^{2}}{2n}
\sum_{p=0}^{\nu(n)-1}\,\sigma_{p}^{2}\biggr)\biggr\}-1\nonumber\\
&= \sum_{k=1}^{+\infty}\sum_{j=1}^{k}\mathbb E \biggl\{W_{j-1}^{n}
\mathbb E \biggl\{\!\!\biggl(\frac{itX_{j}}{\sqrt{n}}-\frac{t^{2}X_{j}^{2}}{2n}+u(\frac{tX_{j}}{\sqrt{n}})
+\frac{t^{2}\sigma_{j-1}^{2}}{2n}+\beta(\!\frac{t^{2}\sigma_{j-1}^{2}}{2n}) \frac{t^{4}\sigma_{j-1}^{4}}{8n^{2}}\biggr)I_{\nu(n)=k}|
\mathcal F_{j-1}\biggr\}\biggr\}
\end{align}
Since $\{\nu(n)=k\}$ is measurable with respect to the $\sigma$-algebra $\mathcal F_{k}$, for all $j\ge 2$, we have
$$\sum_{k=1}^{j-1}\mathbb E \{X_{j}I_{\nu(n)=k}|\mathcal F_{j-1}\}=
\sum_{k=1}^{j-1}\mathbb E \{(X_{j}^{2}-\sigma_{j-1}^{2})I_{\nu(n)=k}|\mathcal F_{j-1}\}=0.$$
On the other hand, for all $j\ge 1$ we have
$$\sum_{k=1}^{+\infty}\mathbb E \{X_{j}I_{\nu(n)=k}|\mathcal F_{j-1}\}=
\sum_{k=1}^{+\infty}\mathbb E \{(X_{j}^{2}-\sigma_{j-1}^{2})I_{\nu(n)=k}|\mathcal F_{j-1}\}=0.$$
It follows that, for all $j\ge 1$
$$\sum_{k\ge j}\mathbb E \{X_{j}I_{\nu(n)=k}|\mathcal F_{j-1}\}=
\sum_{k\ge j}\mathbb E \{(X_{j}^{2}-\sigma_{j-1}^{2})I_{\nu(n)=k}|\mathcal F_{j-1}\}=0.$$
So, from (\ref {dalib}) we derive
\begin{align}\label{h}
&\quad \biggl|\mathbb E \biggl\{\exp\biggl(\frac{itS_{\nu(n)}}{\sqrt{n}}+\frac{t^{2}}{2n}
\sum_{p=0}^{\nu(n)-1}\,\sigma_{p}^{2}\biggr)\biggr\}-1\biggr| \nonumber\\
&=\biggl|\sum_{k=1}^{+\infty}\sum_{j=1}^{k}\mathbb E \biggl\{W_{j-1}^{n}
\mathbb E \biggl\{\biggl(u(\frac{tX_{j}}{\sqrt{n}})
+\beta(\frac{t^{2}\sigma_{j-1}^{2}}{2n}) \frac{t^{4}\sigma_{j-1}^{4}}{8n^{2}}\biggr)I_{\nu(n)=k}|
\mathcal F_{j-1}\biggr\}\biggr\}\biggr| \nonumber\\
&\leq \sum_{k=1}^{+\infty}\sum_{j=1}^{k}\mathbb E \biggl\{\exp\biggl(\frac{t^{2}}{2n}\sum_{p=0}^{j-1}\,
\sigma_{p}^{2}\biggr)
\mathbb E \biggl\{\biggl(\frac{|t|^{3}|X_{j}|^{3}}{6n^{\frac{3}{2}}}
+\frac{t^{4}\sigma_{j-1}^{4}}{8n^{2}}\biggr)I_{\nu(n)=k}|
\mathcal F_{j-1}\biggr\}\biggr\}.
\end{align}
\par For any $j\ge 2$ and any real function $\psi$ such that $\mathbb E (\psi(X_{k}))<\infty$ for any positive $k$, we have
\begin{align}\label{11}
&\sum_{k=1}^{j-1}\mathbb E \biggl\{\exp\biggl(\frac{t^{2}}{2n}\sum_{p=0}^{j-1}\,\sigma_{p}^{2}\biggr)\mathbb E \biggl\{\psi(X_{j})\,I_{\nu(n)=k}|
\mathcal F_{j-1}\biggr\}\biggr\}\nonumber\\
&=\sum_{k=1}^{j-1}\mathbb E \biggl\{\exp\biggl(\frac{t^{2}}{2n}\sum_{p=0}^{j-1}\,\sigma_{p}^{2}\biggr)\mathbb E \biggl\{\psi(X_{j})|
\mathcal F_{j-1}\biggr\}\,I_{\nu(n)=k}\biggr\}.
\end{align}
On the other hand, for all $j\ge 1$, we have
\begin{align}\label{22}
&\,\,\sum_{k=1}^{+\infty}\mathbb E \biggl\{\exp\biggl(\frac{t^{2}}{2n}\sum_{p=0}^{j-1}\,
\sigma_{p}^{2}\biggr)\mathbb E \biggl\{\psi(X_{j})\,I_{\nu(n)=k}|
\mathcal F_{j-1}\biggr\}\biggr\}\nonumber\\
&=\mathbb E \biggl\{\exp\biggl(\frac{t^{2}}{2n}\sum_{p=0}^{j-1}\,\sigma_{p}^{2}\biggr)\,\psi(X_{j})
\biggr\}\nonumber\\
&=\sum_{k=1}^{+\infty}\mathbb E \biggl\{\exp\biggl(\frac{t^{2}}{2n}\sum_{p=0}^{j-1}
\,\sigma_{p}^{2}\biggr)\mathbb E \biggl\{\psi(X_{j})|\mathcal F_{j-1}\biggr\}I_{\nu(n)=k}
\biggr\}.
\end{align}
It follows from (\ref{11}) and (\ref {22}) that
\begin{align}\label{dal}
&\,\,\,\,\sum_{j=1}^{+\infty}\sum_{k \ge j}\mathbb E \biggl\{\exp\biggl(\frac{t^{2}}{2n}\sum_{p=0}^{j-1}\,
\sigma_{p}^{2}\biggr)\mathbb E \biggl\{\psi(X_{j})\,I_{\nu(n)=k}|
\mathcal F_{j-1}\biggr\}\biggr\}\nonumber\\
&=\sum_{j=1}^{+\infty}\sum_{k\ge j}\mathbb E \biggl\{\exp\biggl(\frac{t^{2}}{2n}\sum_{p=0}^{j-1}\,
\sigma_{p}^{2}\biggr)\mathbb E \biggl\{\psi(X_{j})|\mathcal F_{j-1}\biggr\}\,I_{\nu(n)=k}\biggr\}.
\end{align}
\par Applying (\ref{h}) and (\ref{dal}) for $\psi(x)=|x|^{3}$ we deduce that
\begin{align}\label{hh}
&\biggl|\mathbb E \biggl\{\exp\biggl(\frac{itS_{\nu(n)}}{\sqrt{n}}+\frac{t^{2}}{2n}
\sum_{p=0}^{\nu(n)-1}\,\sigma_{p}^{2}\biggr)\biggr\}-1\biggr| \nonumber\\
&\leq \sum_{k=1}^{+\infty}\sum_{j=1}^{k}\mathbb E \biggl\{\exp\biggl(\frac{t^{2}}{2n}\sum_{p=0}^{j-1}\,
\sigma_{p}^{2}\biggr)
\biggl(\mathbb E \biggl\{\frac{|t|^{3}|X_{j}|^{3}}{6n^{\frac{3}{2}}}|
\mathcal F_{j-1}\biggr\}I_{\nu(n)=k}
+\mathbb E \biggl\{\frac{t^{4}\sigma_{j-1}^{4}}{8n^{2}}I_{\nu(n)=k}|
\mathcal F_{j-1}\biggr\}\biggr)\biggr\}\nonumber\\
&\leq \sum_{k=1}^{+\infty}\sum_{j=1}^{k}\mathbb E \biggl\{\exp\biggl(\frac{t^{2}}{2n}\sum_{p=0}^{j-1}\,
\sigma_{p}^{2}\biggr)
\biggl(\frac{|t|^{3}Y_{j-1}\,\sigma_{j-1}^{2}}{6n^{\frac{3}{2}}}I_{\nu(n)=k}
+\frac{t^{4}\sigma_{j-1}^{4}}{8n^{2}}I_{\nu(n)=k}\biggr)\biggr\}
\end{align}
By the H\"older inequality, for all $j\in\mathbb N^{*}$
$$\sigma_{j-1}^{2}=\mathbb E (X_{j}^{2}|\mathcal F_{j-1})\leq \mathbb E (|X_{j}|^{3}|
\mathcal F_{j-1})^{\frac{2}{3}}\leq
Y_{j-1}^{\frac{2}{3}}\,\sigma_{j-1}^{\frac{4}{3}}
\,\,\,\,\,\,\,\,a.s.,$$ whence
\begin{equation}\label{sigma}
\sigma_{j-1}^{2}\le Y_{j-1}^{2}\,\,\,\,\qquad a.s.
\end{equation}
From (\ref {hh}), (\ref{sigma}) and using the fact that $Y_{k}\geq Y_{j-1}\geq 1$ for all $j\leq k$, we deduce that
\begin{align}\label{alah}
&\biggl|\mathbb E \biggl\{\exp\biggl(\frac{itS_{\nu(n)}}{\sqrt{n}}+\frac{t^{2}}{2n}
\sum_{p=0}^{\nu(n)-1}\,\sigma_{p}^{2}\biggr)\biggr\}-1\biggr|\nonumber\\
&\le \biggl(\frac{|t|^{3}}{6n^{\frac{3}{2}}}
+\frac{t^{4}}{8n^{2}}\biggr)\sum_{k=1}^{+\infty}\sum_{j=1}^{k}\mathbb E \biggl\{ Y_{j-1}^{2}\sigma_{j-1}^{2}
\exp\biggl(\frac{t^{2}}{2n}
\sum_{p=0}^{j-1}\,\sigma_{p}^{2}\biggr)I_{\nu(n)=k}\biggr\}\nonumber\\
&\le \biggl(\frac{|t|^{3}}{6n^{\frac{3}{2}}}
+\frac{t^{4}}{8n^{2}}\biggr)\sum_{k=1}^{+\infty} \mathbb E
\biggl\{Y_{k}^{2}\sum_{j=1}^{k}\sigma_{j-1}^{2}
\exp\biggl(\frac{t^{2}}{2n}
\sum_{p=0}^{j-1}\,\sigma_{p}^{2}\biggr)I_{\nu(n)=k}\biggr\}.
\end{align}
\par To bound up the terms appearing in (\ref {alah}), we will use the following elementary lemma. 
\begin{lemme}\label{lemm}
Let $k\ge 1$, then on the event $\{\nu(n)=k\}$ we have
$$\sum_{j=1}^{k}\exp\biggl(\frac{t^{2}}{2n}
\sum_{p=0}^{j-1}\,\sigma_{p}^{2}\biggr)\frac{t^{2}}{2n}\sigma_{j-1}^{2}
\leq \exp(\frac{t^{2}}{2})\biggl(1+\frac{Y_{k}^{2}t^{2}}{n}\biggr).$$
\end{lemme}
{\bf{Proof of Lemma}}.  On the event $\{\nu(n)=k\}$, we have
\begin{align*}
\exp(\frac{t^{2}}{2})&\ge
\exp\biggl(\frac{t^{2}}{2n}\sum_{p=0}^{k-1}\sigma_{p}^{2}\biggr)
-\exp\biggl(\frac{t^{2}}{2n}\sigma_{0}^{2}\biggr)\\
& \ge \sum_{j=1}^{k-1}\exp\biggl(\frac{t^{2}}{2n}\sum_{p=0}^{j-1}\sigma_{p}^{2}\biggr)
\biggl(\exp({\frac{t^{2}\sigma_{j}^{2}}{2n}})-1\biggr)
\end{align*}
Using the inequality,\,\,\,$\exp(x)-1\ge x$  for all $x\ge 0$, one obtains
\begin{align*}
\!\!\!\!\!\!\!\!\!\!\!\!\!\!\!\!\!\!\!\!\!\!\!\!\!\!\!\!\!\!\!\!\!\!\!\exp(\frac{t^{2}}{2})\ge \sum_{j=1}^{k-1}\exp\biggl(\frac{t^{2}}{2n}
\sum_{p=0}^{j-1}\sigma_{p}^{2}\biggr)\frac{t^{2}}{2n}\sigma_{j}^{2}.
\end{align*}
Therefore
\begin{align*}
&\,\,\,\, \sum_{j=1}^{k-1}\exp\biggl(\frac{t^{2}}{2n}
\sum_{p=0}^{j-1}\sigma_{p}^{2}\biggr)\frac{t^{2}}{2n}\sigma_{j-1}^{2}\\
&\leq \sum_{j=1}^{k-1}\exp\biggl(\frac{t^{2}}{2n}
\sum_{p=0}^{j-1}\sigma_{p}^{2}\biggr)\frac{t^{2}}{2n}(\sigma_{j-1}^{2}-\sigma_{j}^{2})
+\exp(\frac{t^{2}}{2})\\
&=\sum_{j=1}^{k-2}\biggl(\exp\biggl(\frac{t^{2}}{2n}
\sum_{p=0}^{j}\sigma_{p}^{2}\biggr)-\exp\biggl(\frac{t^{2}}{2n}
\sum_{p=0}^{j-1}\sigma_{p}^{2}\biggr)\biggr)\frac{t^{2}}{2n}\sigma_{j}^{2}-
\frac{t^{2}}{2n}\exp\biggl(\frac{t^{2}}{2n}
\sum_{p=0}^{k-2}\sigma_{p}^{2}\biggr)\sigma_{k-1}^{2}\\
&\,\,\,\,\,\,\,+\frac{t^{2}}{2n}
\exp(\frac{t^{2}}{2n}\sigma_{0}^{2})\sigma_{0}^{2}+\exp(\frac{t^{2}}{2})\\
&\le \frac{t^{2}}{2n}Y_{k}^{2}\sum_{j=1}^{k-2}\biggl(\exp\biggl(\frac{t^{2}}{2n}
\sum_{p=0}^{j}\sigma_{p}^{2}\biggr)-\exp\biggl(\frac{t^{2}}{2n}
\sum_{p=0}^{j-1}\sigma_{p}^{2}\biggr)\biggr)+\frac{t^{2}}{2n}Y_{k}^{2}
\exp(\frac{t^{2}}{2n}\sigma_{0}^{2})+\exp(\frac{t^{2}}{2})\\
&\le \biggl(1+\frac{t^{2}}{2n}Y_{k}^{2}\biggr)\exp(\frac{t^{2}}{2}).
\end{align*}
We conclude the proof of the lemma by noting that
$\sigma_{k-1}^{2}\le Y_{k}^{2}$ \, and \,\,$\sum\limits_{p=0}^{k-1}\sigma_{p}^{2}\le n$ \,\,a.s..
\par Finally, according to Lemma \ref {lemm}  and the (\ref{alah}) we get
$$\biggl|\mathbb E \biggl\{\exp\biggl(\frac{itS_{\nu(n)}}{\sqrt{n}}+\frac{t^{2}}{2n}
\sum\limits_{p=0}^{\nu(n)-1}\,\sigma_{p}^{2}\biggr)\biggr\}-1\biggr| \leq
a_{n}\exp(\frac{t^{2}}{2})\biggl(\frac{|t|}{3\sqrt{n}}
+\frac{t^{2}}{4n}+\frac{a_{n}|t|^{3}}{3n^{\frac{3}{2}}}+\frac{a_{n}t^{4}}{4n^{2}}\biggr),$$
where $a_{n}=(\mathbb E Y_{\nu(n)}^{4})^{\frac{1}{2}}$.
\subsection{Proof of the Inequality (8)}
Using (\ref{p}) and the inequality $|1-\exp(-x)|\le x$, for all $x\ge 0$ we see that
\begin{align*}
&\quad \,\,\biggl|\mathbb E \biggl\{\exp\biggl(\frac{itS_{\nu(n)}}{\sqrt{n}}+
\frac{t^{2}}{2n}\sum\limits_{p=0}^{\nu(n)-1}\,\sigma_{p}^{2}\biggr)\biggr\}
-\mathbb E \biggl\{\exp\biggl(\frac{itS_{\nu(n)}}{\sqrt{n}}+
 \frac{t^{2}}{2}\biggr)\biggl\}\biggr|\\
 &=\biggl|\mathbb E \biggl\{\exp\biggl(\frac{itS_{\nu(n)}}{\sqrt{n}}+
 \frac{t^{2}}{2}\biggr)\biggl(\exp\biggl(-\frac{t^{2}}{2n}
 \gamma(n) \sigma_{\nu(n)}^{2}\biggr)-1\biggr)\biggl\}\biggr|\\
&\leq \mathbb E \biggl\{\biggl|1-\exp\biggl(-\frac{t^{2}}{2n}
 \gamma(n) \sigma_{\nu(n)}^{2}\biggr)\biggr|\biggr\}\exp(\frac{t^{2}}{2})\\
 &\le \mathbb E \biggl\{\frac{t^{2}}{2n}|\gamma(n)| \sigma_{\nu(n)}^{2}\biggr\}
\exp(\frac{t^{2}}{2})\\
&\leq (\mathbb E Y_{\nu(n)}^{4})^{\frac{1}{2}}\frac{t^{2}}{2n}\exp(\frac{t^{2}}{2}).
\end{align*}
Therefore (8) holds true.
\par From (7) and (8) we conclude that
$$\biggl|\mathbb E \biggl\{\exp(\frac{itS_{\nu(n)}}{\sqrt{n}})\biggr\}-
\exp(-\frac{t^{2}}{2})\biggr|\leq a_{n}\biggl(\frac{|t|}{3\sqrt{n}}+
\frac{3t^{2}}{4n}+\frac{|t|^{3}}{3n^{\frac{3}{2}}} a_{n}+\frac{t^{4}}{4n^{2}}
 a_{n}\biggr).$$
Using Esseen's theorem, we derive $$\sup_{x\in\mathbb{R}}
\biggl|F_{n}(x)- \phi(x)\biggr| \le \frac{ a_{n}}{\pi}\int_{-y}^{y}\biggl(
\frac{1}{3\sqrt{n}}+
\frac{3|t|}{4n}+\frac{t^{2}}{3n^{\frac{3}{2}}}a_{n}+
\frac{ |t|^{3}}{4n^{2}}a_{n}\biggr)dt+\frac{24}{\pi\sqrt{2\pi}y}.$$
Hence
$$\sup_{x\in\mathbb{R}}
\biggl|F_{n}(x)- \phi(x)\biggr| \le \frac{a_{n}}{\pi}\biggl(
\frac{2y}{3\sqrt{n}}+ \frac{3y^{2}}{4n}+\frac{
2y^{3}}{9n^{\frac{3}{2}}}a_{n}+
\frac{y^{4}}{8n^{2}}a_{n}\biggr)+\frac{24}{\pi\sqrt{2\pi}y}.$$
Choosing $y$ in such a way that $y/\sqrt{n}=1/(ya_{n})$, i.e.
 $y=(n/a_{n}^{2})^{\frac{1}{4}}$, we infer that
\begin{align*}
\sup_{x\in\mathbb{R}}
\biggl|F_{n}(x)- \phi(x)\biggr|
&\leq \frac{{a_{n}}^{\frac{1}{2}}}{\pi {n}^{\frac{1}{4}}}\biggl(11+\frac{3}{4n^{\frac{1}{4}}}+\frac{2}{9n^{\frac{1}{2}}}+
\frac{1}{8n^{\frac{3}{4}}}\biggr).
\end{align*}
The proof of the inequality (2) in theorem is complete.
\subsection{Proof of the Inequality (9)}
Observing that the random events $\{\gamma(n)\leq x\}\cap\{\nu(n)=k\}$ and consequently the
random variables
$\sqrt{\gamma(n)}I_{\nu(n)=k}$ are measurable with respect
to $\mathcal F_{k}$, we find that
\begin{align*}
 &\quad\,\, \biggl|\mathbb E \biggl\{\exp\biggl(\frac{itS_{\nu(n)}}{\sqrt{n}}\biggr)\biggr\}
 -\mathbb E \biggl\{\exp\biggl(\frac{itS'_{\nu(n)}}{\sqrt{n}}\biggr)\biggr\}\biggr|\\
 &=\biggl| \sum_{k=0}^{+\infty}\mathbb E \biggl\{\biggl(
 \exp\biggl(\frac{itS_{k}}{\sqrt{n}}\biggr)-
 \exp\biggl(\frac{itS_{k}}{\sqrt{n}}+
 \frac{it\sqrt{\gamma(n)}}{\sqrt{n}}x_{\nu(n)+1}\biggr)\biggr)I_{\nu(n)=k}\biggr\}\biggr|
 \end{align*}
 \begin{align*}
&\leq \sum_{k=0}^{+\infty}\biggl|\mathbb E \biggl\{\,\exp\biggl(\frac{itS_{k}}{\sqrt{n}}\biggr)
\biggl(1-\exp\biggl(\frac{it\sqrt{\gamma(n)}}{\sqrt{n}}X_{\nu(n)+1}\biggr)\biggr)I_{\nu(n)=k}\biggr\}\biggr|\\
&=\sum_{k=0}^{+\infty}\biggl|\mathbb E \biggl\{\,\exp\biggl(\frac{itS_{k}}{\sqrt{n}}\biggr)
\biggl(-\frac{it\sqrt{\gamma(n)}}{\sqrt{n}}\,X_{k+1}
+\frac{t^{2}}{2n}\gamma(n) X_{k+1}^{2}-u(\frac{t\sqrt{\gamma(n)}}{\sqrt{n}}\,X_{k+1})\biggr)I_{\nu(n)=k}\biggr\}\biggr|\\
&= \sum_{k=0}^{+\infty}\biggl|\mathbb E \biggl\{\exp\biggl(\frac{itS_{k}}{\sqrt{n}}\biggr)
\mathbb E \biggr\{-\frac{it\sqrt{\gamma(n)}}{\sqrt{n}}X_{k+1}+\frac{t^{2}}{2n}\gamma(n) X_{k+1}^{2}
-u(\frac{t\sqrt{\gamma(n)}X_{k+1}}{\sqrt{n}})|\mathcal F_{k}\biggr\}I_{\nu(n)=k}\biggr\}\biggr|\\
&=\sum_{k=0}^{+\infty}\biggl|\mathbb E \biggl\{\exp\biggl(\frac{itS_{k}}{\sqrt{n}}\biggr)
\biggl(\frac{t^{2}}{2n}\gamma(n) X_{k+1}^{2}-u(\frac{t}{\sqrt{n}}\sqrt{\gamma(n)}
\,X_{k+1})\biggr)I_{\nu(n)=k}\biggr\}\biggr|\\
&\le \sum_{k=0}^{+\infty}\mathbb E \biggl\{I_{\nu(n)=k}\,\frac{3t^{2}}{2n}\gamma(n) X_{k+1}^{2}\biggr\}\\
&\leq \frac{3t^{2}}{2n} \sum_{k=0}^{+\infty}\mathbb E \biggl\{I_{\nu(n)=k}\,\,\mathbb E \{ X_{k+1}^{2}|\mathcal F_{k}\}\biggr\}\\
&\le \frac{3t^{2}}{2n}\, \mathbb E (Y_{\nu(n)}^{4})^{\frac{1}{2}}.
\end{align*}
The proof of the inequality (9) is complete.
\subsection{Proof of the inequality (3).}
According to Esseen's theorem where $y=(n/a_{n}^{2})^{\frac{1}{4}}$
\,\,and the inequality (9), one obtains
\begin{align*}
\sup_{x\in\mathbb {R}}\biggl|H_{n}(x)-\phi(x)\biggr|
&\leq \frac{1}{\pi}
\int_{-y}^{y}\biggl|\mathbb E \biggl\{\exp(\frac{itS'_{\nu(n)}}{\sqrt{n}})\biggr\}-
\exp(-\frac{t^{2}}{2})\biggr|\,\frac{dt}{|t|}\,+ \frac{24}{\pi \sqrt{2\pi}y}\\
&\le \frac{{a_{n}}^{\frac{1}{2}}}{\pi {n}^{\frac{1}{4}}}\biggl(11+\frac{3}{4n^{\frac{1}{4}}}+\frac{2}{9n^{\frac{1}{2}}}+
\frac{1}{8n^{\frac{3}{4}}}\biggr)+\frac{1}{\pi}\int_{-y}^{y}
\frac{3|t|}{2n}E(Y_{\nu(n)}^{4})^{\frac{1}{2}}\,dt\\
&\leq \frac{{a_{n}}^{\frac{1}{2}}}{\pi {n}^{\frac{1}{4}}}\biggl(11+\frac{3}{4n^{\frac{1}{4}}}+\frac{2}{9n^{\frac{1}{2}}}+
\frac{1}{8n^{\frac{3}{4}}}\biggr)+\frac{3}{2\pi\sqrt{n}}\\
&\leq \frac{{a_{n}}^{\frac{1}{2}}}{\pi {n}^{\frac{1}{4}}}\biggl(11+\frac{9}{4n^{\frac{1}{4}}}+\frac{2}{9n^{\frac{1}{2}}}+
\frac{1}{8n^{\frac{3}{4}}}\biggr).
\end{align*}
The proof of theorem is complete.   \quad\quad $\square$\\
\par Proofs of corollaries 1, 2 and 3 are easy so, it is left to the reader.\\
\par {\textbf{Acknowledgements:}} The author thanks the referee
for careful reading of the manuscript and for valuable suggestions
which improved the presentation of this paper.

\end{document}